\documentclass[english,twoside,11pt]{article}

\usepackage[german,english]{babel}
\usepackage{amsmath}
\usepackage{amssymb,latexsym}
\usepackage{url}
\usepackage{graphicx}
\usepackage{booktabs}
\usepackage{psfrag}
\usepackage{paralist}

\newtheorem{deff}{Definition}

\newtheorem{prop}[deff]{Proposition}
\newtheorem{thm}[deff]{Theorem}

\newtheorem{conj}[deff]{Conjecture}

\makeatletter

\title{Surface Realization with the Intersection Segment Functional}

\author{\Large Stefan Hougardy\footnote{Supported by the DFG Research Center {\sc Matheon} ``Mathematics for Key Technologies", Berlin}, 
               Frank H.~Lutz\footnote{Supported by the DFG Research Group ``Polyhedral Surfaces'', Berlin}, 
               and Mariano Zelke$^*$}

\date{}

\begin{document}

\selectlanguage{english}

\maketitle

\begin{abstract}
Deciding realizability of a given polyhedral map on a (compact, connected) surface
belongs to the hard problems in discrete geometry, from the theoretical,
the algorithmic, and the practical point of view. 

In this paper, we present a heuristic algorithm for the realization of simplicial maps, 
based on the intersection segment functional. The heuristic was used
to find geometric realizations in ${\mathbb R}^3$ for \emph{all} 
vertex-minimal triangulations of the orientable surfaces
of genus $g=3$ and $g=4$. Moreover, for the first time, examples of simplicial polyhedra 
in ${\mathbb R}^3$ of genus $5$ with $12$ vertices were obtained. 
\end{abstract}

\section{Introduction}

A \emph{polyhedral map} on a surface is a (finite) set of polygons
(with at least three sides), which are glued together (topologically) 
along edges to form the surface, such that there are no self-identifications 
on the boundaries of the polygons, and two polygons are either disjoint
or intersect in exactly one edge or one vertex only. 
We thus can think of a polyhedral map as a combinatorial model for a surface.

For a given polyhedral map it is natural to try to visualize it as 
a \emph{polyhedron} in three-space or in higherdimensional space ${\mathbb R}^d$,
such that every polygon is the convex hull of its vertices
and two polygons are either disjoint in ${\mathbb R}^d$,
they intersect in a common edge and are not coplanar, 
or they intersect in a common vertex only. Such a realization
usually is called a \emph{geometric} or \emph{polyhedral realization},
with straight edges, plane polygons, and no non-trivial intersections
(with neighboring polygons being not coplanar).

\medskip

\begin{figure}
\begin{center}
\small
\psfrag{1}{1}
\psfrag{2}{2}
\psfrag{3}{3}
\psfrag{4}{4}
\psfrag{5}{5}
\psfrag{6}{6}
\psfrag{7}{7}
\psfrag{8}{8}
\includegraphics[width=.7\linewidth]{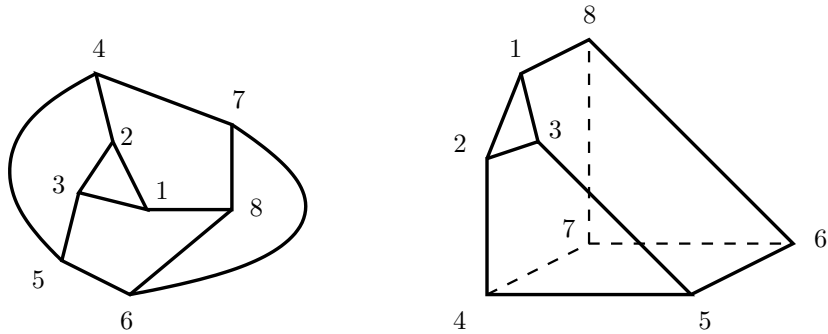} 
\end{center}
\caption{A polyhedral map on $S^2$ and a corresponding geometric realization in ${\mathbb R}^3$.}
\label{fig:map}
\end{figure}

\noindent
\emph{Example 1:} A polyhedral map on the $2$-sphere $S^2$ consisting 
of the polygons $123$, $12478$, $13568$, $2354$, $4567$, and $678$ 
together with a corresponding realization in ${\mathbb R}^3$ 
is displayed in Figure~\ref{fig:map}.

\medskip

Realizability of maps on the $2$-sphere $S^2$
was proved by Steinitz (\cite{Steinitz1922},
    \cite{SteinitzRademacher1934}; cf.\ also
    \cite[Ch.~13]{Gruenbaum1967}, \cite[Lec.~4]{Ziegler1995}):
Every polyhedral map on the $2$-sphere $S^2$
is geometrically realizable in~${\mathbb R}^3$ 
as the boundary complex of a convex $3$-polytope.

However, not all polyhedral maps are realizable.
For example, \emph{simple polyhedral maps} (i.e., maps with all vertices of valence three)
on surfaces different from the $2$-sphere $S^2$ are not realizable in any\, ${\mathbb R}^d$
(see Gr\"unbaum \cite[Ex.\ 11.1.7, Ex.\ 13.2.3]{Gruenbaum1967}).

\medskip

\begin{figure}
\begin{center}
\small
\psfrag{1}{1}
\psfrag{2}{2}
\psfrag{3}{3}
\psfrag{4}{4}
\psfrag{5}{5}
\psfrag{6}{6}
\psfrag{7}{7}
\psfrag{8}{8}
\psfrag{9}{9}
\psfrag{10}{10}
\psfrag{11}{11}
\psfrag{12}{12}
\psfrag{13}{13}
\psfrag{14}{14}
\includegraphics[width=.9\linewidth]{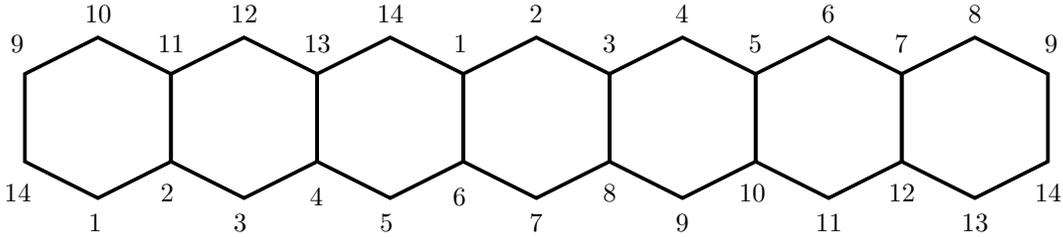} 
\end{center}
\caption{The non-realizable $6$-$3$-equivelar map with $14$ vertices on the torus.}
\label{fig:moebius_dual}
\end{figure}

\noindent
\emph{Example 2:} All \emph{$6$-$3$-equivelar} maps on the torus 
(i.e., maps consisting of only $6$-gons with every vertex
belonging to exactly three $6$-gons) are simple and therefore cannot be realized
in any ${\mathbb R}^d$.
The smallest example (see Figure~\ref{fig:moebius_dual})
of the family is the combinatorial dual of M\"obius' $7$-vertex triangulation
of the torus \cite{Moebius1886}. A ``realization'' of this $6$-$3$-torus
with flat, but non-convex $6$-gons was given by Szilassi
\cite{Szilassi1986}, the \emph{Szilassi-torus}.

\medskip

Betke and Gritzmann \cite{BetkeGritzmann1982} found a further combinatorial obstruction to geometric
realizability: Let $W$ be any subset of the set of odd valent vertices of a
polyhedral map $M^2$ and let $F_W$ be the set of facets containing some vertex of $W$.
If $2|F_W|\leq |W|$, then $M^2$ is not realizable in any ${\mathbb R}^d$.
Again, the Betke-Gritzmann obstruction rules out realizability of $6$-$3$-equivelar maps on the torus,
but the obstruction was also used by McMullen, Schulz, and Wills \cite{McMullenSchulzWills1982} 
to show non-realizability for other, non-simple families of equivelar maps.

In a \emph{simplicial map} (i.e., a triangulation of a surface as a simplicial complex)
every triangle contains at most three odd valent vertices,
from which it can be deduced that $|F_W|\geq |W|$ for every
subset $W$ of odd valent vertices. Thus, the Betke-Gritzmann obstruction 
cannot be applied to simplicial maps to show non-realizability. 

Gr\"unbaum \cite{Gruenbaum1967} proposed in 1967 that every triangulated torus  
should be geometrically realizable in ${\mathbb R}^3$. His famous conjecture
was open for 40 years and was proved in 2007 by Archdeacon, Bonnington, and
Ellis-Monaghan \cite{ArchdeaconBonningtonEllisMonaghan2007}.
For the class of triangulations of the projective plane with one face
removed geometric realizability in ${\mathbb R}^3$ was proved in 2008 
by Bonnington and Nakamoto \cite{BonningtonNakamoto2008}. 
However, not every triangulated M\"obius strip needs to be 
realizable in ${\mathbb R}^3$, a~counterexample is due to by Brehm \cite{Brehm1983}.

Until rather recently, no computational tools were available
\emph{to actually find realizations} for simplicial surfaces.
In the past 30 years the most promising approach to obtain a polyhedral realization in ${\mathbb R}^3$
for a given triangulation was to try to \emph{build a physical model}, 
for example, by exploiting symmetries or by employing the \emph{rubber band technique} 
of Bokowski~\cite{Bokowski2008}.

Very basic heuristic procedures for finding realizations in ${\mathbb R}^3$ for larger
classes of examples were used by Fendrich \cite{Fendrich2003}
(to show that all triangulated tori with up to 11 vertices are
realizable via embeddings in the $2$-skeleta of random $4$-polytopes)
and Lutz \cite{Lutz2008a} (to obtain realizations for triangulations 
of the orientable surface of genus $2$ via choosing coordinates randomly).
These methods, although useful in processing larger numbers of examples,
are less powerful than human imagination leading to a physical model.
For example, $864$ of the $865$ vertex-minimal triangulations 
of the orientable surface of genus $2$ were realized with 
the random realization method in \cite{Lutz2008a}, 
with a total computation time of 30 CPU months.
However, it was not possible to randomly realize the remaining example 
(with the highest combinatorial symmetry of the 865 examples). 
The computer search was run unsuccessfully for more than a month 
before the example was realized within a day by Bokowski \cite{Bokowski2008}
with the rubber band method.

In this paper, we present a heuristic algorithm for finding polyhedral realizations 
for (closed, orientable) simplicial surfaces in ${\mathbb R}^3$, which, for the first time, 
surpasses the physical approach with respect to its processing time 
and its qualitative range of examples. In particular, we show that 
all vertex-minimal triangulations of the orientable surfaces 
of genus $g=3$ and $g=4$ are realizable. We also provide examples 
of vertex-minimal simplicial polyhedra of genus $5$ with $12$ vertices.

In the following section we give a brief survey 
on realizability results for surfaces, 
vertex-minimal triangulations, and algorithmic aspects of deciding realizability. 
Section~\ref{sec:functional} is devoted to our realization heuristic 
based on the \emph{intersection segment functional}.
Computational results are presented in Section~\ref{sec:computational}.
An extension of our approach to convex realizations of triangulated
spheres is discussed in Section~\ref{sec:convex}.

\section{Realizability of Polyhedral Surfaces \\ and Polyhedral Complexes}
\label{sec:realizability}

In general, every $d$-dimensional simplicial complex (with $n$ vertices) is 
polyhedrally embeddable in ${\mathbb R}^{2d+1}$, as it can be
realized as a subcomplex of the boundary complex of the cyclic
polytope $C(n,2d+2)$; cf.\ Gr\"unbaum \cite[Ex.\ 25, p.\ 67]{Gruenbaum1967}.
However, van Kampen \cite{vanKampen1932} and Flores~\cite{Flores1933} showed that
$d$-dimensional simplicial complexes 
cannot always be embedded topologically in~${\mathbb R}^{2d}$,
e.g., the $d$-skeleton ${\rm Sk}_d(\Delta_{2d+2})$
of the $(2d+2)$-simplex $\Delta_{2d+2}$ is not embeddable
in~${\mathbb R}^{2d}$. For further examples and references
see Matou\v{s}ek \cite[5.1]{Matousek2003}, Novik~\cite{Novik2000}, 
and Schild~\cite{Schild1993}.

For smooth $d$-manifolds, Whitney \cite{Whitney1944}
proved that they can smoothly be embedded in~${\mathbb R}^{2d}$,
and Penrose, Whitehead, and Zeeman \cite{PenroseWhiteheadZeeman1961} 
showed that for $0<2(k+1)\leq d$
every $k$-connected PL (i.e., piecewise linear) 
$d$-manifold has a PL embedding in  ${\mathbb R}^{2d-k}$.
In particular, surfaces have PL embeddings in ${\mathbb R}^4$.
Orientable surfaces (with or without boundary) and non-orientable
surfaces with boundary are even PL embeddable in~${\mathbb R}^3$
(which follows from the classification of surfaces by Dehn and Heegaard \cite{DehnHeegaard1907}).
Closed non-orientable $d$-manifolds cannot be embedded 
topologically in ${\mathbb R}^{d+1}$; cf.\ Bredon \cite[p.\ 353]{Bredon1993}.

Thus, for triangulated orientable surfaces 
(with or without boundary) and for triangulated 
non-orientable surfaces with boundary we have:
\begin{compactitem}
\item PL embeddability in ${\mathbb R}^3$
\item and polyhedral realizability in ${\mathbb R}^5$.
\end{compactitem}
Triangulations of closed non-orientable surfaces are
\begin{compactitem}
\item not (topologically) embeddable in ${\mathbb R}^3$,
\item but are PL embeddable in ${\mathbb R}^4$,
\item and are polyhedrally realizable in ${\mathbb R}^5$.
\end{compactitem}

Perles showed (cf.~\cite[11.1.8]{Gruenbaum1967})
that a polyhedral map is realizable in some ${\mathbb R}^d$
if and only if it is realizable in ${\mathbb R}^5$.

A natural approach to establish geometric realizability in ${\mathbb R}^3$
for polyhedral maps on orientable surfaces of genus $g\geq 1$ 
is to identify a given polyhedral map as a subcomplex 
of the boundary complex of a convex $4$-polytope $P$.
The Schlegel diagram of $P$ then yields coordinates
for the realization in ${\mathbb R}^3$; see, for example, 
McMullen, Schulz, and Wills \cite{McMullenSchulzWills1982} 
for realizations of equivelar maps obtained this way,
and cf.\ Altshuler \cite{Altshuler1971b,Altshuler1971} 
for combinatorial properties on maps that guarantee 
realizability via Schlegel diagrams.
 
Altshuler and Brehm \cite{AltshulerBrehm1984} gave
a polyhedral map $T_8$ on the torus with only $8$ vertices
which is realizable in ${\mathbb R}^3$ (cf.\ also Simutis \cite{Simutis1977}), 
but \emph{not} via the Schlegel diagram of a convex $4$-polytope. In fact, 
the map $T_8$ is not isomorphic to a subcomplex of the boundary complex 
of any convex polytope~\cite{AltshulerBrehm1984}.

Realizability (via subcomplexes of convex $5$-polytopes) 
of triangulations of the torus and the projective plane in ${\mathbb R}^4$ 
was proved by Brehm and Schild \cite{BrehmSchild1995}, herewith
sharpening Barnette's result \cite{Barnette1983} on the geometric realizability
of triangulations of the projective plane in~${\mathbb R}^4$.

Polyhedral surfaces that are obtained by projections (of $2$-dimensional 
subcomplexes) of higherdimensional polytopes together with obstructions 
to projectability are discussed by R\"orig, Sanyal, 
and Ziegler~\cite{Roerig2009,RoerigSanyal2009pre,RoerigZiegler2009pre}.
Knotted realizations of triangulated tori are studied by Lutz,
Sullivan, and Witte in \cite{LutzWitte2007pre}.
For further results and references on polyhedral maps see Brehm and Wills \cite{BrehmWills1993},
Brehm and Schulte \cite{BrehmSchulte1997}, and Ziegler \cite{Ziegler2008b}.

\subsection{Simplicial Maps}

Let $M$ be a (closed) triangulated surface with $n=f_0$ vertices, $f_1$ edges,
and $f_2$ triangles, i.e., $M$ has \emph{face-vector} $f=(n,f_1,f_2)$.
If $M$ has Euler characteristic $\chi(M)$, then by Euler's equation,
$$n-f_1+f_2 = \chi(M).$$
Double counting of the incidences between edges and triangles of the triangulation 
yields $2f_1=3f_2$. So together,
$$f=(n,3n-3\chi(M),2n-2\chi(M)).$$
A triangulated surface with $n$ vertices obviously has at most $f_1\leq \binom{n}{2}$ edges.
By plugging in $f_1=3n-3\chi(M)$ we obtain
Heawood's bound~\cite{Heawood1890} from 1890
that a triangulation of a $2$-manifold $M$ of Euler characteristic $\chi(M)$
has at least
\begin{equation}
n\geq\Bigl\lceil\tfrac{1}{2}(7+\sqrt{49-24\chi (M)})\Bigl\rceil
\end{equation}
vertices. Heawood's bound is sharp for all surfaces, except for
the orientable surface of genus~$2$, the Klein bottle, 
and the non-orientable surface of genus~$3$, where an extra vertex 
has to be added to the lower bound.

Corresponding \emph{vertex-minimal triangulations} (i.e., triangulations 
with the minimal possible number of vertices)
of the real projective plane ${\mathbb R}{\bf P}^2$ with $6$ vertices 
and of the $2$-torus with $7$ vertices (M\"obius' torus \cite{Moebius1886})
were already known in the 19th century,
but it took until 1955 to complete the construction of series
of examples of vertex-minimal triangulations for all non-orientable surfaces 
(Ringel~\cite{Ringel1955}) and until 1980 for all orientable surfaces 
(Jungerman and Ringel~\cite{JungermanRingel1980}).

If a given triangulation of an orientable surface
is realizable in ${\mathbb R}^3$, then so are subdivisions
of it that are obtained by successively subdividing edges
and triangles. Hence, vertex-minimal triangulations apparently
are good candidates for non-realizable simplicial maps.
Hereby, triangulations with 
\begin{equation}\label{eq:neighborly}
n=\tfrac{1}{2}(7+\sqrt{49-24\chi (M)})
\end{equation}
are of particular interest (cf.\ \cite{Csaszar1949}), as for these we have $f_1=\binom{n}{2}$,
that is, the respective triangulations are \emph{neighborly} with
complete $1$-skeleton (which should make realizability difficult).

A polyhedral realization of the combinatorially unique 
vertex-minimal $7$-vertex triangulation of the torus with $f=(7,21,14)$
was given by Cs\'asz\'ar~\cite{Csaszar1949}, \cite{Lutz2002b}
(although realizability possibly was known already to M\"obius;
cf.\ \cite[p. 553]{Moebius1886}, \cite{Reinhardt1885}).

The next case of equality in (\ref{eq:neighborly}) 
yields $59$ examples of vertex-minimal $12$-vertex triangulations 
of the orientable surface of genus $6$
\cite{AltshulerBokowskiSchuchert1996}; see below.

\subsection{Realizability vs.\ Non-Realizability of Simplicial Maps}

For every individual triangulation of an orientable surface,
realizability (in ${\mathbb R}^3$) can be decided algorithmically by
the following two-step
procedure, cf.\ \cite{Bokowski2001}, \cite[Ch.\ VIII]{BokowskiSturmfels1989}:
\begin{compactitem}
\item[1.] Enumerate all oriented matroids compatible with the given triangulation.
      If there are none, then the triangulation is not realizable, else
\item[2.] decide realizability of the oriented matroids from 1.\ via
      solving associated polynomial inequality systems. 
\end{compactitem}
Theoretically, the second step can be done algorithmically
(for example, with Collins' \emph{Cylindrical Algebraic Decomposition}
algorithm \cite{Collins1975}). In practice, however, there are no methods known 
that would work sufficiently fast to yield results even for small examples.
See \cite{Bokowski2001} and \cite[Ch.\ VIII]{BokowskiSturmfels1989} for more 
comments on this and on algebraic tools such as final polynomials.
For general polyhedral maps on orientable surfaces Brehm proved that
the realizability problem is NP-hard (as a consequence of his universality 
theorem for realization spaces of maps, cf.~\cite{Ziegler2008a}). 
The complexity of the realization problem restricted to simplicial maps 
is unknown. In fact, it was open for a long time, 
whether there are non-realizable examples at all.

In a major breakthrough, Bokowski and Guedes de Oliveira \cite{BokowskiGuedes_de_Oliveira2000}
showed in 2000 (using 10 CPU years) that one of the $59$ vertex-minimal 
$12$-vertex triangulations of the orientable surface of genus~$6$ 
has no compatible orientable matroid and therefore is not realizable.

Schewe \cite{Schewe2007,Schewe2008pre} substantially improved the enumeration
of compatible orientable matroids and was able to verify that,
in fact, all $59$ vertex-minimal $12$-vertex triangulations 
of the orientable surface of genus $6$ are non-realizable.
Moreover, he found three examples of non-realizable vertex-minimal
$12$-vertex triangulations of the orientable surface of genus~$5$.
At~least for one of these examples it is possible to remove a
triangle from the triangulation while maintaining non-realizability.
Connected sums with other triangulations then still are non-realizable.
Hence, for every orientable surface of genus $g\geq 5$ there are triangulations 
that cannot be realized geometrically in ${\mathbb R}^3$.

Apart from the approach via oriented matroids, non-realizability results 
(for simplicial maps in ${\mathbb R}^3$)
seem to be difficult to achieve: Novik \cite{Novik2000} 
associated an integer program with a given triangulation, which, 
if it has no solution, yields non-realizability. 
Improved systems have been proposed by Timmreck \cite{Timmreck2008}. 
So far, however, all tested systems for orientable surfaces 
either had solutions or turned out to be computationally intractable. 
In a different approach, Brehm \cite{Brehm1983} used a linking number argument 
to show that there is a non-realizable triangulation of the M\"obius strip
with $9$ vertices.

\subsection{Heuristics for the Realization of Simplicial Maps}

Until recently, it was considered to be rather difficult and time-consuming 
to actually find realizations for given triangulations. 
Examples of polyhedral realizations of vertex-minimal triangulations
of the orientable surfaces of genus $3$ and $4$ with $10$ and $11$ 
vertices, respectively, were constructed \emph{by hand} by Brehm \cite{Brehm1981,Brehm1987b}
and Bokowski and Brehm \cite{BokowskiBrehm1987,BokowskiBrehm1989}. 
Some of these examples were found by exploiting combinatorial symmetries
of the triangulations, others with the \emph{rubber band technique} of Bokowski \cite{Bokowski2008}.

A simple computer heuristic (\emph{by choosing coordinates randomly})
was used in~\cite{Lutz2008a} to show that $864$ of the $865$ examples of vertex-minimal
triangulations of the orientable surface of genus $2$ are realizable.
The remaining case then was settled by Bokowski with the rubber band method \cite{Bokowski2008}.
All $865$ examples were later found to have realizations with 
\emph{small coordinates}~\cite{HougardyLutzZelke2007a}, 
i.e., all these examples are realizable with integer coordinates in general
position in the $(4\times 4\times 4)$-cube. Moreover, 
realizations in the $(5\times 5\times 5)$-cube were obtained 
for $17$ of the $20$ vertex-minimal triangulations with $10$ vertices 
of the orientable surface of genus $3$ \emph{by isomorphism-free enumeration} 
of possible coordinate configurations in general position~\cite{HougardyLutzZelke2007b}.

In the following, we will discuss an improved heuristic to obtain
polyhedral realizations in ${\mathbb R}^3$ for triangulations of 
orientable surfaces. In particular, we will show that all
vertex-minimal triangulations of the orientable surfaces of genus $g=3$ and $g=4$
are realizable and that there are examples of simplicial polyhedra of genus $5$
with $12$ vertices.

\section{Realization with the Intersection Segment Functional}
\label{sec:functional}

As mentioned in the previous section there have been so far three major
heuristics for the realization of simplicial surfaces 
(of genus $g\geq 1$) in ${\mathbb R}^3$:
\begin{compactitem}
\item by explicit geometric construction
      \cite{BokowskiBrehm1987,BokowskiBrehm1989,Brehm1981,Brehm1987b} 
      (e.g., via the rubber band technique of Bokowski~\cite{Bokowski2008});
\item by choosing coordinates randomly \cite{Lutz2008a};
\item by enumeration of realizations with small coordinates
     \cite{HougardyLutzZelke2007a,HougardyLutzZelke2007b,HougardyLutzZelke2008}.
\end{compactitem}

As a more sophisticated approach we suggest to proceed as follows.
For a given triangulation (of an orientable surface 
of small genus with few vertices)

\begin{compactitem}
\item[1.] start with random coordinates for the vertices of the
          triangulation
\item[2.] and then ``move vertices around'' to eventually obtain a
          realization.
\end{compactitem}

\pagebreak

\begin{figure}
\begin{center}
\small
\psfrag{a}{$a$}
\psfrag{b}{$b$}
\psfrag{c}{$c$}
\psfrag{d}{$d$}
\psfrag{e}{$e$}
\psfrag{f}{$f$}
\psfrag{r}{$r$}
\psfrag{s}{$s$}
\psfrag{p}{$p$}
\psfrag{q}{$q$}
\includegraphics[height=50mm]{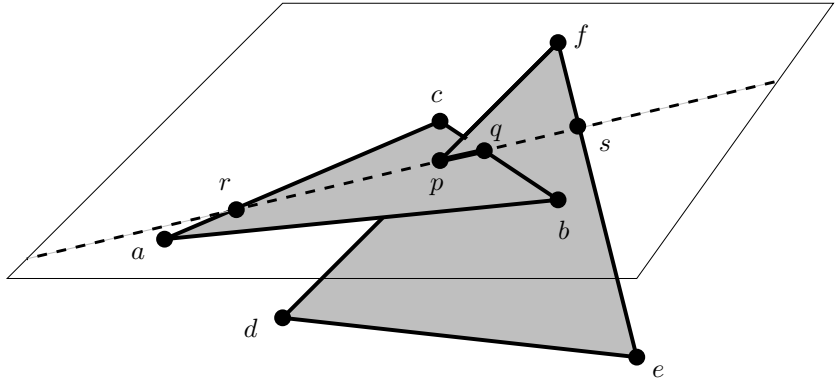}
\end{center}
\caption{Two intersecting triangles.}
\label{fig:two_triangles}
\end{figure}

\noindent
For the second step we take as an objective to minimize
the \emph{intersection segment functional}:
\begin{quote}
\emph{Let $M^2$ be a triangulated orientable surface with vertex-set $V$ and let
$V_{{\mathbb Z}^3}$ be a set of\, $|V|$ integer vertices in general position
in ${\mathbb R}^3$. Then every pair of triangles of $M^2$ coordinatized
with the coordinates of $V_{{\mathbb Z}^3}$ either has empty
intersection in ${\mathbb R}^3$ or intersects in a segment;
see Figure~\ref{fig:two_triangles} for the intersection segment $p$--$q$ of two triangles.
The sum of the lengths of the intersection segments over all pairs of (non-neighboring) triangles 
is the \emph{intersection segment functional}.
}
\end{quote}
We require that the points are \emph{in general position},
i.e., no three points are on a line and no four points are on a plane,
in order to avoid degenerate intersections of triangles.
Further, we use integer coordinates and therefore move the points
in the second step above on the integer grid ${\mathbb Z}^3$ only. 
Our aim then will be to find integer coordinates in general position
for which the intersection segment functional vanishes
for the given triangulation. 

From an initial set of random coordinates we proceed to minimize 
the intersection segment functional by a local search of \emph{hill-climbing} type:
\begin{quote}
\emph{In every step, we randomly pick a vertex $v\in V_{{\mathbb Z}^3}$ 
and a coordinate direction, $\pm x$, $\pm y$, or $\pm z$, 
and then move the vertex $v$ one integer step into the respective direction.
If the resulting set of coordinates is in general position
and the new value of the intersection segment functional is strictly smaller
than before, the move is accepted and the next step is executed. 
Otherwise the move is discarded and we start anew from the previous set of coordinates.}

\emph{If all possible choices of moves have been tested for some set of coordinates
without improvement, then we are stuck in a \emph{local minimum}.
In this case, for one step only, we choose one of the \emph{admissible moves}, 
i.e., a move that yields a set of coordinates in general position,
but which not necessarily decreases the intersection segment functional.
From there, we then try to continue to decrease the intersection segment
functional in a new direction.}
\end{quote}

\pagebreak

\begin{figure}
\begin{center}
\includegraphics[height=51mm]{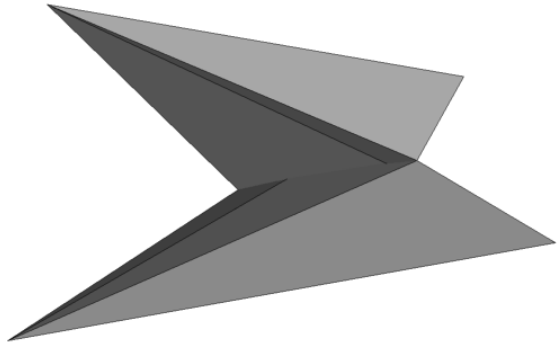} 
\end{center}
\caption{A ``non-realization'' of the octahedron with locally minimal functional.}
\label{sphereLocalMin}
\end{figure}

\noindent
\emph{Example 3:} Local minima with positive value of the intersection segment
functional can occur even for small triangulations. For example, 
the boundary of the octahedron with triangles
{\small
\begin{center}
\begin{tabular}{@{\extracolsep{3.5mm}}llllllll}
$123$ & $124$ & $135$ & $145$ & $236$ & $246$ & $356$ & $456$ 
\end{tabular} 
\end{center}
}
\noindent
and furnished with coordinates (in general position)
{\small
\begin{center}
\begin{tabular}{l@{\hspace{2mm}}l@{\hspace{5mm}}l@{\hspace{2mm}}l@{\hspace{5mm}}l@{\hspace{2mm}}l@{\hspace{5mm}}l@{\hspace{2mm}}l@{\hspace{5mm}}l@{\hspace{2mm}}l@{\hspace{5mm}}l@{\hspace{2mm}}l@{\hspace{5mm}}l@{\hspace{2mm}}l}
1: & (4,4,6)  &  2: & (5,6,6)  &  3: & (9,7,4)  &  4: & (5,9,1)  &  5:  &  (4,6,3)  &  6:  &  (1,5,7)
\end{tabular} 
\end{center}
}
\noindent
attains a local minimum for the intersection segment functional with value 3.17;
see Figure~\ref{sphereLocalMin} for a visualization.

\subsection{Details of the Algorithm}

Initially the vertices of the triangulation are placed randomly at
general positions in a cube of size $50\times50\times50$.
This cube is chosen at the center of a larger $(250\times250\times250)$-cube 
that we take as \emph{bounding box} for all possible positions of the
vertices during the local search. 
After the choice of the starting positions the smaller cube is not used anymore. 
\begin{compactitem}
\item Thus, we allow the diameter of the vertex-set to increase moderately (which possibly 
      helps to decrease the intersection segment functional by unfolding the initial shape).
\item At the same time there is a fixed lower bound for the change, at every step, of the intersection segment functional
      (determined by the size of the bounding box and the fact that we admit integer coordinates only).
      This way we avoid that the sequence of improvements for the functional converges to zero.
\end{compactitem}
An \emph{admissible step} then is a movement of one vertex by one integer in one
of the coordinate directions such that the resulting set of coordinates is in general position
and is within the bounding box.

\begin{compactitem}
\item If the intersection segment functional becomes zero, a realization for the given triangulation is found. 
\item If a realization is not found within a fixed period of time $T$,
      the whole process is restarted for the triangulation, beginning with the random selection 
      of the starting coordinates (in the smaller cube). In doing so we try to overcome situations 
      in which the process cycles between different local minima.
\end{compactitem}
A standard problem with local search algorithms is to appropriately choose the
parameters that govern the procedure. For some of the 20 examples of vertex-minimal
$10$-vertex triangulations of the orientable surface of genus $3$ 
we tried the following variants:
\begin{compactitem}
\item We chose different sizes for the initial cube, ranging from $5\times5\times5$ to
      $500\times500\times500$. 
\item We allowed the bounding box to be between 1 up to 8 times the size of the initial cube.
\item If the segment functional decreases by moving one vertex in one direction, 
      we moved the vertex as far as possible in this direction 
      (until the intersection segment functional starts to increase again).
\item In case of a local minimum we determined all pairs of vertices
      for which the exchange of their positions decreases the
      intersection segment functional. We then executed one such exchange
      at random. If there is no such pair, we randomly exchanged
      two arbitrary vertices.
\item Instead of minimizing the intersection segment functional we tried 
      to minimize the \emph{normalized intersection segment functional},
      which is obtained from the intersection segment functional
      by dividing by the total length of the edges of the coordinatized triangulation.
\item We first generated 10000 sets of initial coordinates of which
      we selected the set with the smallest functional before starting the
      local search.
\end{compactitem}
From all these variants the previously described one turned out to have the best performance.
This variant then was used to find realizations for other triangulations.

\subsection{Test Sets of Minimal Triangulations}

\begin{figure}
\begin{center}
\includegraphics[height=50mm]{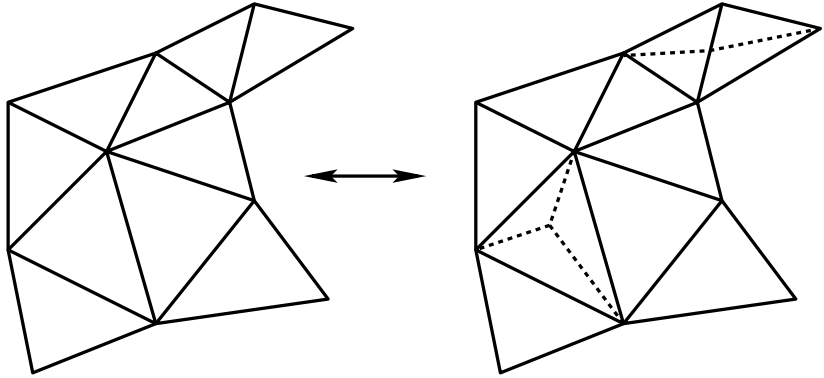}
\end{center}
\caption{Subdivisions of a triangle and of an edge.}
\label{fig:subdivision}
\end{figure}

If some triangulation of an orientable surface is realizable, then so are all subdivisions of it
that result from the starting triangulation by an iterative sequence of \emph{elementary subdivisions} 
of triangles and edges (see Figure~\ref{fig:subdivision}).
For every geometrically realized triangulation in such a sequence, we can always place the new vertex 
slightly ``above'' or ``below'' the respective triangle or the respective edge
of the previous realization. Alternatively, we could choose all new vertices 
on the original surface and then slightly perturb the coordinates of
the new vertices into general position.

\pagebreak

\begin{figure}
\begin{center}
\includegraphics[height=35mm]{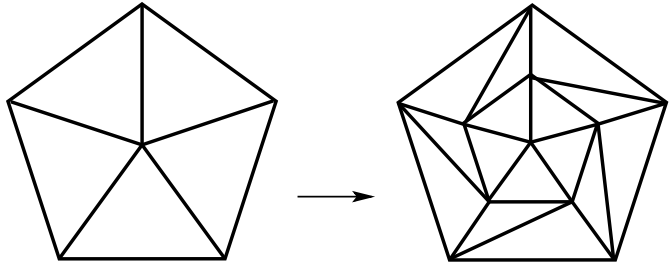}
\end{center}
\caption{Replacement of the star of a vertex of degree $5$.}
\label{fig:pentagon_replace}
\end{figure}

A triangulation of a surface is \emph{minimal} if it does not result 
from a triangulation with fewer vertices by a sequence of elementary 
subdivisions. If all minimal triangulations of a surface are realizable,
then all triangulations of the surface are realizable.
Unfortunately, for surfaces of genus $g\geq 1$ the set of minimal
triangulations is infinite: it comprises the infinite set of triangulations 
with all vertices of degree at least~$5$, since for any such triangulation
we can replace the \emph{star} of a vertex (i.e., all triangles that
contain the vertex) with a patch that has more vertices, 
but all of degree at least~$5$; see Figure~\ref{fig:pentagon_replace}. 
Equivelar triangulations of surfaces of genus $g\geq 1$ are minimal.
For $g=1$ there are infinitely many equivelar triangulations,
whereas for $g>1$ there are only finitely many examples; cf.~\cite{SulankeLutz2006pre}.

A finite subset (of the set of minimal triangulations) of particular interest 
to test realizability is the set of vertex-minimal triangulations. 
If these are realizable, then this should give a strong indication that, 
in fact, all triangulations of the surface are realizable.

\begin{figure}
\begin{center}
\includegraphics[height=52mm]{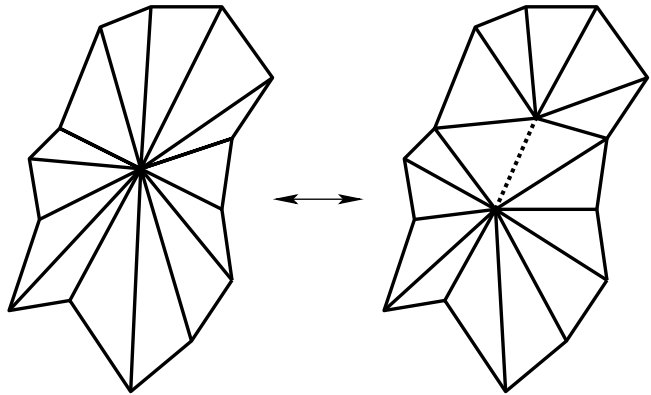}
\end{center}
\caption{Expansion (respectively contraction) of an edge.}
\label{fig:edge_contraction}
\end{figure}

A larger, but still finite set of minimal triangulations that contains all vertex-minimal triangulations
is defined as follows: If we allow \emph{edge expansions} (with \emph{edge contractions} 
as inverses; see Figure~\ref{fig:edge_contraction}) instead of elementary subdivisions, 
then for every surface there is only a finite set (see Barnette and Edelson~\cite{BarnetteEdelson1988}) 
of \emph{irreducible} triangulations for which no edge can be contracted
without changing the topological type of the triangulation.
Unfortunately, it is, a priori, not clear whether a realizable triangulation of an orientable surface 
of genus $g\geq 1$ remains realizable after the expansion of an edge. 
At least, for every explicit polyhedral realization it can easily be tested 
whether a particular edge expansion can be carried out (via a system 
of linear constraints on the link of the respective vertex; we thank
the anonymous referee for pointing this out to us).

It follows from the work of Steinitz \cite[\S 46]{SteinitzRademacher1934}
that every triangulated $2$-sphere can be reduced to the boundary
of the tetrahedron by a sequence of edge contractions, that is,
the boundary of the tetrahedron is the only irreducible
triangulation of the $2$-sphere.
Gr\"unbaum and Lavrenchenko \cite{Lavrenchenko1990} determined
the  number of irreducible triangulations of the torus: there are $21$ such examples 
with up to $10$ vertices and they are all realizable.
Sulanke~\cite{Sulanke2006apre,Sulanke2006bpre,Sulanke2005cpre} showed by enumeration that there
are exactly $396.784$ examples of irreducible triangulations (with up to $17$ vertices) 
of the orientable surface of genus~$2$.

Although it might be desirable to test realizability for a larger set
of irreducible triangulations, we restricted ourselves to vertex-minimal ones. 
There is only one unique vertex-minimal triangulation of the torus, i.e.,
M\"obius $7$-vertex torus \cite{Moebius1886} for which Cs\'asz\'ar~\cite{Csaszar1949}
gave an explicit polyhedral model. Vertex-minimal triangulations of the orientable surfaces
of genus $2$ and~$3$ were enumerated in \cite{Lutz2008a}, those of genus $4$ and $5$ in 
\cite{SulankeLutz2006pre}, and the vertex-minimal examples of genus $6$ in  \cite{AltshulerBokowskiSchuchert1996};
see Table~\ref{tbl:vertex_minimal} for the corresponding minimal numbers of vertices $n_{\rm min}$
and the respective numbers of combinatorial types of triangulations.

\begin{table}
\caption{Numbers of vertex-minimal triangulations of the orientable surfaces of genus $g\leq 6$.}\label{tbl:vertex_minimal}
\centering
\defaultaddspace=0.15em
\begin{tabular}{@{}r@{\hspace{8mm}}r@{\hspace{8mm}}r@{}}
\\\toprule
 \addlinespace
 \addlinespace
 \addlinespace
 \addlinespace
  $g$  &   $n_{\rm min}$  &  Types \\
\midrule
\\[-4mm]
 \addlinespace
 \addlinespace
 \addlinespace
 \addlinespace
 
 \addlinespace
0 &  4 &       1   \\
 \addlinespace
1 &  7 &       1   \\
 \addlinespace
2 & 10 &     865   \\
 \addlinespace
3 & 10 &      20   \\
 \addlinespace
4 & 11 &     821   \\
 \addlinespace
5 & 12 & 751.593   \\
 \addlinespace
6 & 12 &      59   \\
 \addlinespace

 \addlinespace
 \addlinespace
 \addlinespace
 \addlinespace
\bottomrule
\end{tabular}
\end{table}

\section{Computational Results}
\label{sec:computational}

\subsubsection*{\rm\em Genus $2$}
In~\cite{Lutz2008a}, geometric realizations for $864$ of the $865$ vertex-minimal $10$-vertex 
triangulations of the orientable surface of genus~$2$ were found with the random realization approach 
in a total computation time of $30$~CPU months on a 2.8 GHz processor; the remaining example then 
was realized with the rubber band method~\cite{Bokowski2008}. For realizations of the $865$ examples 
with small coordinates see~\cite{HougardyLutzZelke2007a} and the comments above. With our new heuristic
algorithm, based on the intersection segment functional, realizations for the $865$ triangulations were obtained 
in a total time of 218 CPU minutes on a 3.5 GHz processor.

\subsubsection*{\rm\em Genus $3$}
Realizations for $5$ of the $20$ vertex-minimal $10$-vertex triangulations of the orientable 
surface of genus $3$ were constructed by hand by Brehm and Bokowski \cite{BokowskiBrehm1987,Brehm1981,Brehm1987b}.
The random realization approach of~\cite{Lutz2008a}, however, produced no results for these $5$ 
(and for the other~$15$) examples, where we stopped the search after one CPU week each. 
Therefore, the basic random realization approach is not suitable for triangulations of
surfaces of higher genus (or with more vertices). For $17$ of the $20$ triangulations realizations 
with small coordinates in the $(5\times 5\times 5)$-cube were obtained in \cite{HougardyLutzZelke2007b};
this search was run (in total) for 2 CPU years on a 3.5 GHz processor.
Thus, the first task for our new program was to realize the remaining three examples.

\begin{thm}
All $20$ vertex-minimal $10$-vertex triangulations of the orientable surface of \linebreak
genus~$3$
are geometrically realizable in ${\mathbb R}^3$.
\end{thm}
Sets of coordinates for the realizations are available online at \cite{HougardyLutzZelke2007b,Lutz_PAGE}.
In total, it took 28~CPU hours on a 3.5 GHz processor to realize the 20 examples with the help of the intersection
segment functional. For two of the last three of the $20$ examples, we later found realizations
in the $(6\times 6\times 6)$-cube; cf.\ \cite{Lutz_PAGE}.

\subsubsection*{\rm\em Genus $4$}
A first example of a polyhedron of genus $4$ with $11$ vertices was described by Bokowski and Brehm \cite{BokowskiBrehm1989}.
With our intersection segment functional algorithm we found realizations for $626$ 
of the $821$ vertex-minimal $11$-vertex triangulations of the orientable surface of genus $4$. 
In an effort to speed up the search, realizations for the remaining $195$ triangulations 
were obtained \emph{by recycling of coordinates}, that is, whenever a new realizations 
was found we tried to reuse the respective set of coordinates for other triangulations. 
We also slightly distorted the coordinates and then tried to use these coordinates 
for other triangulations; see \cite{Lutz2008a} for additional comments.

\begin{figure}
\begin{center}
\psfrag{13 }{13}\psfrag{15 }{15}\psfrag{17 }{17}\psfrag{19 }{19}
\psfrag{21 }{21}\psfrag{23 }{23}\psfrag{25 }{25}
\psfrag{50 }{50}\psfrag{75 }{75}\psfrag{100 }{100}
\psfrag{125 }{125}\psfrag{150 }{150}
\psfrag{natural logarithm of the number of used local search steps}{natural logarithm 
        of the number of used local search steps}
\psfrag{frequency of found realizations}{frequency of found realizations}
\includegraphics[width=0.65\textwidth]{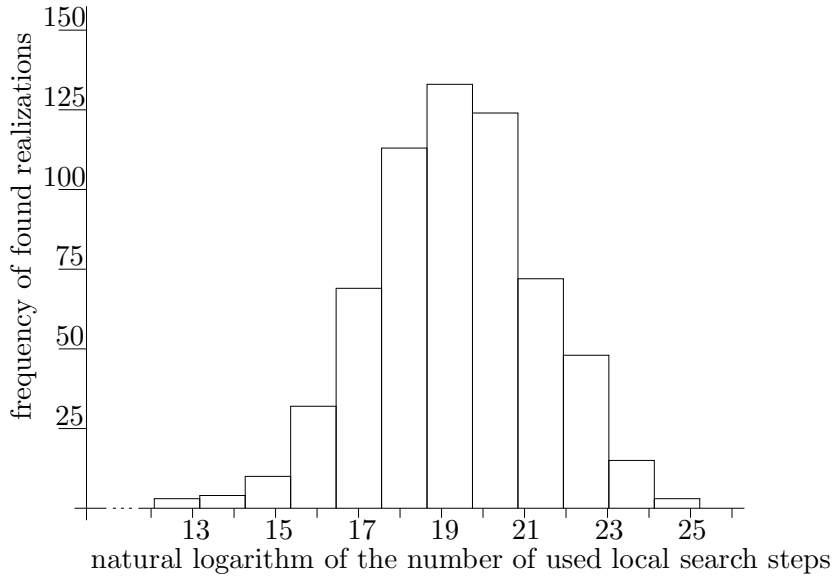}
\end{center}
\caption{Histogram of the natural logarithms of the used local search steps for $626$ realizations of genus $4$.}
\label{stepHistogram}
\end{figure}

\begin{thm}
All $821$ vertex-minimal $11$-vertex triangulations of the orientable surface of genus $4$
are geometrically realizable in ${\mathbb R}^3$.
\end{thm}
We needed a total of $9.51\cdot10^{11}$ steps of the local search process
to realize the $626$ triangulations. As time interval $T$ we chose $15$ minutes, 
so if after $15$ minutes (about $5.4\cdot10^{6}$ steps) a~realization was not reached, 
the search was restarted with new initial coordinates.

Figure \ref{stepHistogram} displays a histogram of the natural logarithms of the number 
of used steps. The picture indicates that the logarithms of the used steps 
are normally distributed, i.e., the used steps underlie a log-normal distribution. 
To confirm this, we ran as a goodness-of-fit-test \cite[Ch.~30]{StuartOrd1991}
the Anderson-Darling test (cf. \cite[p.~10]{Everitt1998}). 
The test estimates the mean to be $19.3$ and the standard deviation to be $2$. 
It yields a $p$-value of $0.5$, which is far above the rejection value of $0.05$. 
Therefore we can view the logarithms of the used steps 
to be normally distributed with the estimated parameters.

Our implementation of the local search process is performing
about $3.6\cdot10^{5}$ steps per minute on a 3.5 GHz processor. 
Therefore, we needed a total of 5 CPU years to realize all triangulations. 
On average, it took 2.9 CPU days for finding a realization for a single triangulation. 

\subsubsection*{\rm\em Genus $5$}
As mentioned in Section~\ref{sec:realizability}, Schewe \cite{Schewe2007,Schewe2008pre} showed
that there are at least three examples of vertex-minimal $12$-vertex triangulations
of the orientable surface of genus $5$ that cannot be realized geometrically in ${\mathbb R}^3$.

\begin{figure}
\begin{center}
\includegraphics[height=61mm]{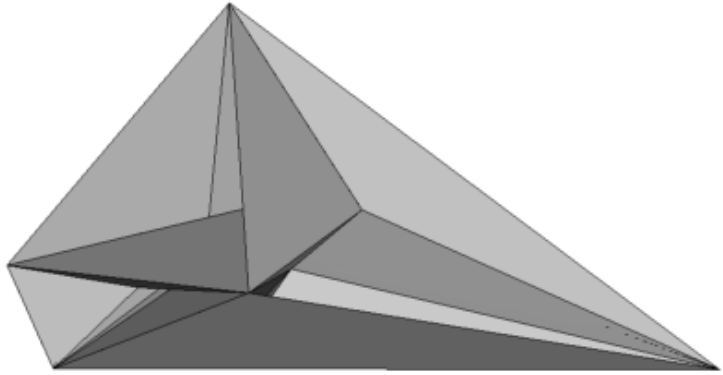}
\end{center}
\caption{A polyhedron of genus $5$.}
\label{fig:genus_5}
\end{figure}

In order to complement Schewe's result, we tried to find realizations 
for at least some of the $751.593$ triangulations.
To this aim we started our process on randomly selected triangulations out of all the $751.593$ vertex-minimal triangulations. 
If after $15$ minutes a realization was not found, a new triangulation was selected at random. 
This way, we tried about 94.000 triangulations, using a total of
$7.52\cdot10^{11}$ local search steps 
-- a CPU time of approximately 4 years -- and succeeded in realizing 15 triangulations. 

\begin{thm}
At least $15$ of the $751.593$ vertex-minimal $12$-vertex triangulations of the orientable surface of genus $5$
are geometrically realizable in ${\mathbb R}^3$.
\end{thm}
Since the 94.000 triangulations we tested were chosen randomly from the list of the $751.593$
genus $5$ triangulations, probably at least $120$ (and perhaps many more) of the examples are realizable.

\bigskip

\noindent
\emph{Example 4:} Figure~\ref{fig:genus_5} displays one of the
polyhedra of genus $5$ with $12$ vertices, which has triangles
{\small
\begin{center}
\begin{tabular}{@{\extracolsep{2.75mm}}llllllllll}
$123$    & $124$    & $135$    & $146$    & $157$       & $168$    & $179$       & $18\,10$    & $19\,10$    & $236$ \\
$245$    & $258$    & $26\,10$ & $28\,11$ & $29\,11$    & $29\,12$ & $2\,10\,12$ & $35\,11$    & $368$       & $378$ \\
$37\,10$ & $39\,10$ & $39\,11$ & $459$    & $46\,11$    & $478$    & $47\,12$    & $489$       & $4\,10\,11$ & $4\,10\,12$ \\
$569$    & $56\,10$ & $57\,10$ & $58\,12$ & $5\,11\,12$ & $679$    & $67\,12$    & $6\,11\,12$ & $89\,12$    & $8\,10\,11$ 
\end{tabular}
\end{center}
}
\noindent
and coordinates
{\small
\begin{center}
\begin{tabular}{l@{\hspace{2mm}}l@{\hspace{5mm}}l@{\hspace{2mm}}l@{\hspace{5mm}}l@{\hspace{2mm}}l@{\hspace{5mm}}l@{\hspace{2mm}}l@{\hspace{5mm}}l@{\hspace{2mm}}l}
1: & (137,124,141) &  2: & (107,118,143) &  3: & (132,130,125) &  4: & (122,127,129) \\
5: & (124,129,132) &  6: & (126,130,124) &  7: & (126,129,129) &  8: & (122,125,138) \\ 
9: & (124,128,136) & 10: & (119,133,134) & 11: & (120,130,135) & 12: & (121,128,133).
\end{tabular} 
\end{center}
}
\noindent
The coordinates for the other $14$ examples can be found online at~\cite{Lutz_PAGE}.

\medskip

Combining the result of Schewe \cite{Schewe2007,Schewe2008pre} (that there are non-realizable triangulations
of the orientable surface of genus $5$) with our finding (that all vertex-minimal
triangulations of surfaces of genus $g\leq 4$ are realizable) gives rise to:

\begin{conj}
Every triangulation of an orientable surface of genus
$g\leq 4$ is geometrically realizable.
\end{conj}
The conjecture holds for genus $0$ (Steinitz \cite{Steinitz1922,SteinitzRademacher1934})
and for genus $1$ (Archdeacon, Bonnington, and Ellis-Monaghan~\cite{ArchdeaconBonningtonEllisMonaghan2007}).

\subsection{Examples with More Vertices}

We also tried our program on some triangulations of tori with more vertices. 
It turned out that it still is possible to find realizations, although 
it takes much longer for every step of the local search process: 
There are $O(|V|^2)$ pairs of triangles that have to be considered 
for the computation, respectively for the update, of the 
intersection segment functional. Moreover, there are $6|V|$ 
possible moves from a current set of coordinates that lead to 
a new set of coordinates. In the worst case, we are forced 
to test almost all these moves just to carry out a single improvement step.
Finally, the initial value of the intersection segment functional
will be larger for triangulations with more vertices,
thus, forcing us to perform more steps.

\pagebreak

\begin{figure}
\begin{center}
\includegraphics[height=33mm]{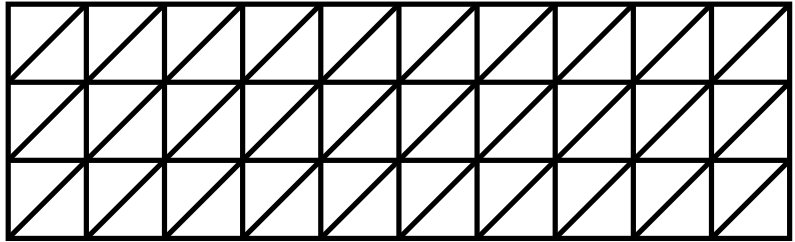}\\[5mm]
\end{center}
\caption{The standard $(3\times 10)$-torus.}
\label{fig:standard_2d}
\end{figure}

\begin{figure}
\begin{center}
\includegraphics[height=67mm]{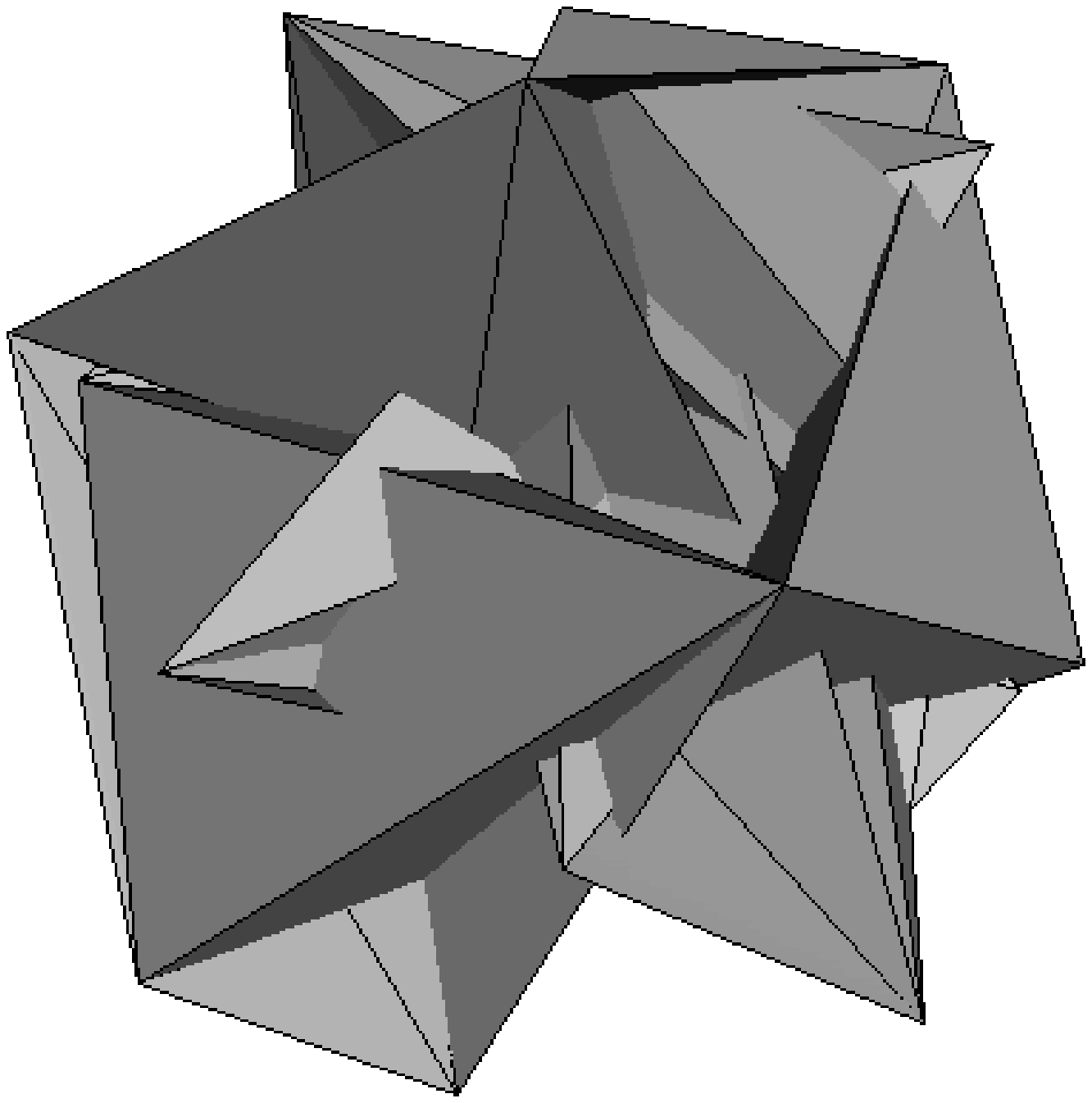}\hspace{5mm}\includegraphics[height=67mm]{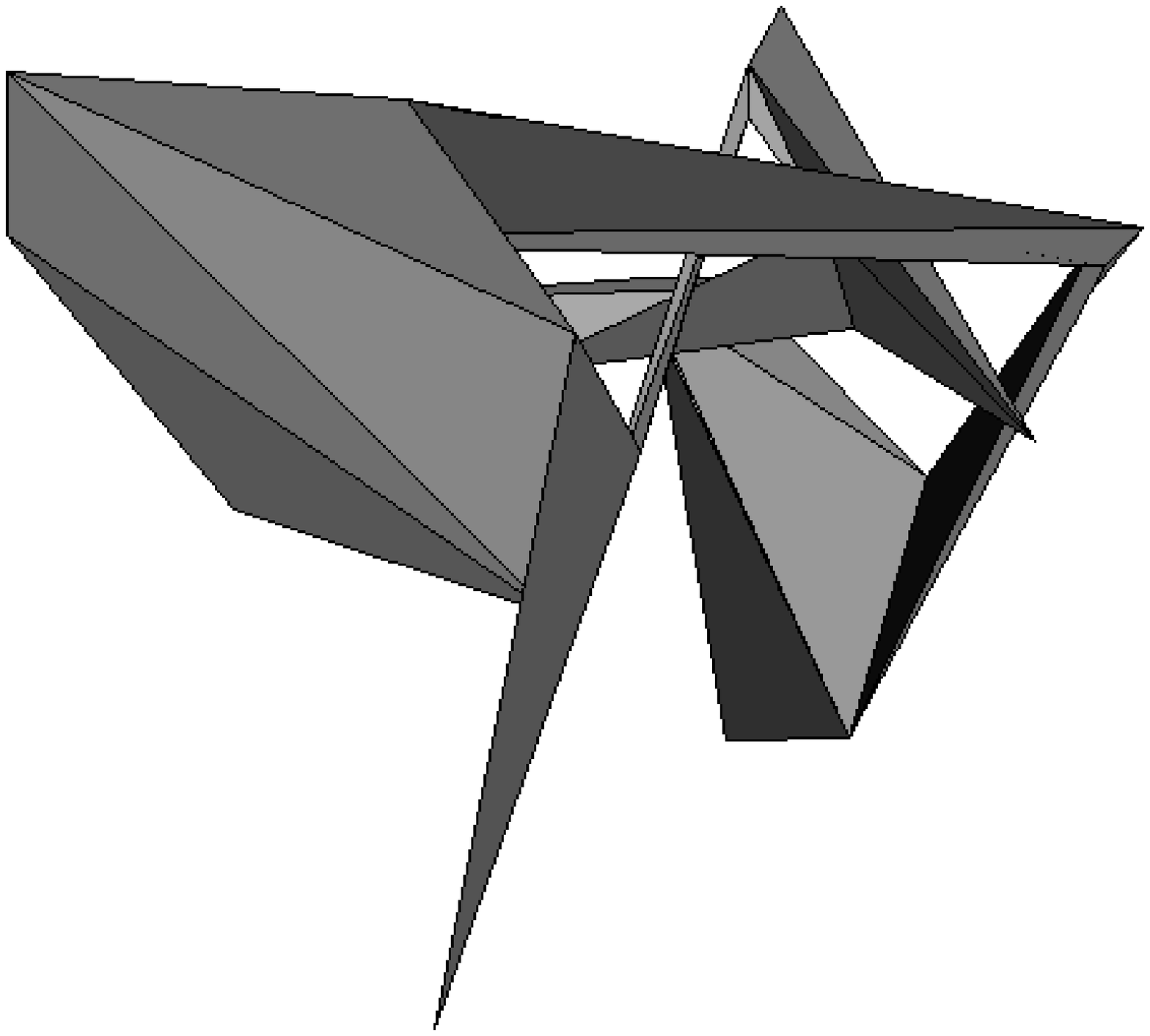}
\end{center}
\caption{The standard $(3\times 10)$-torus with random coordinates and a proper realization.}
\label{fig:standard_3d}
\end{figure}

\noindent
\emph{Example 5:} For the standard $(3\times 10)$-torus (Figure~\ref{fig:standard_2d})
we started with random coordinates (Figure~\ref{fig:standard_3d}, left)
and an initial value $7924.26$ of the intersection segment functional.
It then took 3042 local search steps to obtain a proper realization
(Figure~\ref{fig:standard_3d}, right).

\section{Convex Realizations of Triangulated 2-Spheres}
\label{sec:convex}

According to Steinitz \cite{Steinitz1922,SteinitzRademacher1934}, 
every polyhedral map on the $2$-sphere $S^2$ is geometrically 
realizable in~${\mathbb R}^3$ as the boundary complex 
of a convex $3$-polytope. Tutte's equilibrium method~\cite{Tutte1963}
(see also \cite[Section 12.2]{Richter_Gebert1996}) allows to obtain
corresponding realizations algorithmically via first constructing
a planar equilibrium representation of the edge graph of a given map.
The resulting planar graph can then be interpreted as the Schlegel diagram
of a $3$-polytope. 

Triangulated $2$-spheres are realizable as boundary complexes of simplicial 
$3$-polytopes. However, simplicial $3$- and higher-dimensional spheres 
need not be polytopal. The Br\"uckner-Gr\"unbaum $3$-sphere~\cite{GruenbaumSreedharan1967} 
and the Barnette $3$-sphere \cite{Barnette1973c}, both with $8$-vertices, 
are the smallest examples of non-polytopal simplicial spheres.

In the following, we give a simple modification of our realization heuristic 
in order to obtain \emph{convex realizations} of triangulated $2$-spheres.
(By using an \emph{intersection area functional} one might generalize
this approach to search for convex realizations of simplicial $3$-spheres
in ${\mathbb R}^4$.)

A \emph{non-face} of a triangulated $2$-sphere with $n$ vertices
is a two-element subset (an edge) or a three-element subset (a triangle) 
of the ground set of $n$ vertices that does not constitute a face 
of the triangulation. In any convex realization (in general position)
of the triangulated $2$-sphere as the boundary complex of a simplicial 
$3$-polytope, every non-face intersects with the interior of the 
respective polytope. In particular, every face of the boundary $2$-sphere 
either has no intersection with a given non-face or intersects the non-face 
in an edge or a vertex of the boundary sphere.

Hence, by adding to the intersection segment functional the lengths
of intersection segments for all pairs of triangles consisting
of a triangle of the triangulation and a triangle that does
not belong to the triangulation 
the resulting \emph{extended intersection segment functional}
can be used to obtain convex realizations for triangulated $2$-spheres.
To be more precise:

\begin{prop}
If a triangulated $2$-sphere has no vertex of degree $3$,
then the extended intersection segment functional is zero 
if and only if a convex realization (with vertices 
in general positions) has been reached.
\end{prop}

\noindent
\textbf{Proof.} If a convex realization has been reached,
then obviously the extended intersection segment functional
is zero.

For the other direction, assume that the functional is zero 
and that the vertices are in general position.
In case some vertex $v$ is contained in the convex hull
of the other $n-1$ vertices, then $v$ is contained
in the convex hull of some subset $\{v_1,v_2,v_3,v_4\}$ 
of four of the $n-1$ vertices: Pick any vertex $v_1$ on the 
boundary of the convex hull of the $n$ vertices, then by the 
general position assumption, there is a unique triangle $\{v_2,v_3,v_4\}$
(opposite to $v_1$ with respect to $v$) on the boundary of the convex hull 
such that the tetrahedron $\{v_1,v_2,v_3,v_4\}$ contains $v$.
Without loss of generality, we may assume that no other of the $n-1$ vertices 
is contained in the tetrahedron $\{v_1,v_2,v_3,v_4\}$:
If there is an additional such vertex, say, $v'$, then the convex hull of $v'$ with the 
triangle on the boundary of the tetrahedron $\{v_1,v_2,v_3,v_4\}$ opposite to $v'$
with respect to $v$ is a tetrahedron (of smaller volume) that contains~$v$.
Moreover, since the smallest example of a triangulated $2$-sphere without a vertex of degree $3$
is the boundary complex of the octahedron with $6$ vertices, there is at 
least one vertex of the triangulation that lies outside the tetrahedron $\{v_1,v_2,v_3,v_4\}$.
If the vertex $v$ has degree larger than $4$, then at least one of the triangles 
of the star of $v$ intersects non-trivially some boundary triangle (a face or a non-face 
of the triangulation) of the tetrahedron $\{v_1,v_2,v_3,v_4\}$. 
The line segment the two triangles intersect in contributes a positive value
to the extended intersection segment functional, contradiction.
If $v$ has degree~$4$, then there are two cases. If the star of the vertex $v$ 
contains a vertex different from $v_1,v_2,v_3,v_4$, then again at least one of 
the triangles in the star of $v$ intersects some triangle of the boundary 
of the tetrahedron $\{v_1,v_2,v_3,v_4\}$, contradiction.
Else, let $v_5$ be a vertex outside the convex hull of the vertices $v_1,v_2,v_3,v_4$.
Then $v$ lies in the convex hull of $v_5$ and some triangle  $\{v_{i_1},v_{i_2},v_{i_3}\}$
of the tetrahedron $\{v_1,v_2,v_3,v_4\}$. But then the vertex star of $v$ contains the vertex $v_4$, 
which lies outside the tetrahedron spanned by the vertices $v_{i_1},v_{i_2},v_{i_3},v_5$. 
This again leads to a contradiction.\hfill $\Box$

\pagebreak

In case a triangulation has vertices of degree $3$, non-convex realizations 
of the triangulation can have vanishing extended intersection segment functional.
The smallest such example is the bipyramid over a triangle with one apex pushed inside 
the convex hull of the other four vertices.

Nevertheless, we can recursively remove vertices of degree $3$ from a given triangulation. 
The resulting triangulation then either is the boundary of a tetrahedron or a triangulation 
with vertices all of degree at least four. After obtaining a realization 
for the simplified triangulation, the removed vertices can be added back 
by placing them suitably ``above'' the triangles which they subdivide. 

We successfully tested our approach for some small triangulations of $S^2$: 
there are $233$ triangulations of $S^2$ with $10$ vertices of which $12$ examples
have no vertices of degree $3$. It took, on average, about $5$ minutes 
to obtain convex realizations for these $12$ triangulations.

\medskip

\noindent
\emph{Remark:} Although our main focus in this paper was on the realization of closed, orientable
               triangulated surfaces, the intersection segment functional can, of course, also be used 
               to search for realizations in three-space for other $2$-dimensional
               simplicial complexes. Moreover, the functional can easily be modified for a
               search for immersions of (orientable or non-orientable) triangulated surfaces.

\subsection*{Acknowledgment}
We are grateful to the anonymous referee for many helpful comments.

\bibliography{}

\vspace{10mm}

\small

\noindent
Stefan Hougardy \\
Universit\"at Bonn\\
Forschungsinstitut f\"ur Diskrete Mathematik\\
Lenn\'estr.\ 2\\
53113 Bonn\\
Germany\\
{\tt hougardy@or.uni-bonn.de}

\vspace{10mm}

\noindent
Frank H. Lutz\\
Technische Universit\"at Berlin\\
Institut f\"ur Mathematik \\
Stra{\ss}e des 17.\ Juni 136\\
10623 Berlin\\
Germany\\
{\tt lutz@math.tu-berlin.de}

\vspace{10mm}

\noindent
Mariano Zelke\\
Humboldt-Universit\"at zu Berlin\\
Institut f\"ur Informatik\\
Unter den Linden 6\\
10099 Berlin\\
Germany\\
{\tt zelke@informatik.hu-berlin.de}

\end{document}